\documentclass[10pt,reqno]{elsarticle}

 \usepackage{graphicx,amssymb,amstext,amsmath}
 \usepackage{enumerate}
 \usepackage{amsthm}
\usepackage{xcolor}
\usepackage[colorlinks=true,
linkcolor=blue,
urlcolor=blue,
citecolor=blue]{hyperref}
\usepackage[margin = 2.3cm]{geometry}
\usepackage{listings}\usepackage{float}

\newtheoremstyle{plainsl}%
	{\topsep}
	{\topsep}
	{\slshape} % only non-default setting
	{}
	{\normalfont\bfseries}
	{.}
	{ }
	{}
	
\theoremstyle{plainsl}	
 \newtheorem{theorem}{Theorem}[section]
 \newtheorem{thm}[theorem]{Theorem}

 \newtheorem{lem}[theorem]{Lemma}
 
 \newtheorem{prop}[theorem]{Proposition}
 \newtheorem{quest}[theorem]{Question}
 \newtheorem{conj}[theorem]{Conjecture}
 
 \newtheorem{example}[theorem]{Example}

\DeclareMathOperator\Aut{Aut}

\DeclareMathOperator\sym{Sym}
\DeclareMathOperator\alt{Alt}

\DeclareMathOperator\Alt{Alt}

\DeclareMathOperator\fld{\mathbb{F}}

\DeclareMathOperator\elsm{sum}

\DeclareMathOperator\tr{tr}

\DeclareMathOperator\der{der}
\DeclareMathOperator\id{id}
\newcommand{\agl}[2]{\operatorname{AGL}_#1(#2)}
\newcommand{\pgl}[2]{\operatorname{PGL}_#1(#2)}
\newcommand{\psl}[2]{\operatorname{PSL}_#1(#2)}
\newcommand{\asl}[2]{\operatorname{ASL}_#1(#2)}

\newcommand{\agammal}[2]{\operatorname{A\Gamma L}_#1(#2)}
\newcommand{\pgammal}[2]{\operatorname{P\Gamma L}_#1(#2)}

\newcommand{\mathieu}[1]{\operatorname{M}_{#1}}

\begin{document}

\title{On the Intersection Density of the Kneser Graph $K(n,3)$}

\author[Uregina]{Karen Meagher\fnref{fn1}}
\ead{karen.meagher@uregina.ca}
 
\author[address1,address2]{Andriaherimanana Sarobidy Razafimahatratra\corref{cor1}}
\cortext[cor1]{Corresponding author}

\address[Uregina]{Department of Mathematics and Statistics, University of Regina, Regina, Saskatchewan S4S 0A2, Canada} 

\address[address1]{University of Primorska, UP FAMNIT, Glagolja\v{s}ka 8, 6000 Koper, Slovenia\\
       %University of Primorska, UP IAM, Muzejski trg 2, 6000 Koper, Slovenia 
       }\ead{sarobidy@phystech.edu}
\address[address2]{University of Primorska, UP IAM, Muzejski trg 2, 6000 Koper, Slovenia }

\fntext[fn1]{Research supported by NSERC Discovery Research Grant, Application No.: RGPIN-2018-03952.}

	\date{\today}

\begin{abstract}
  A set $\mathcal{F} \subset \sym(V)$ is \textsl{intersecting} if any two of its elements agree on some element of $V$. Given a finite transitive permutation group $G\leq \sym(V)$, the \textsl{intersection density} $\rho(G)$ is the maximum ratio $\frac{|\mathcal{F}||V|}{|G|}$ where $\mathcal{F}$ runs through all intersecting sets of $G$. The \textsl{intersection density} $\rho(X)$ of a vertex-transitive graph $X = (V,E)$ is equal to $\max \left\{ \rho(G) : G \leq \Aut(X), \mbox{ $G$ transitive} \right\}$. In this paper, we study the intersection density of the Kneser graph $K(n,3)$, for $n\geq 7$.  
  The intersection density of $K(n,3)$ is determined  whenever its automorphism group contains $\psl{2}{q}$, with some exceptional cases depending on the congruence of $q$.  We also briefly consider the intersection density of $K(n,2)$ for values of $n$ where $\psl{2}{q}$ is a subgroup of its automorphism group.
\end{abstract}
\maketitle

%%%%%%%%%%%%%%%%%%%%%%%%%%%%%%%%%%%%%%
%%%%%%%%%%%%%%%%%%%%%%%%%%%%%%%%%%%%%%

\section{Introduction}

There have been many recent papers looking at the size of the largest set of intersecting permutations in a transitive permutation group, see for example~\cite{ ellis2011intersecting, HKKMMO, HKMMwreath, HMMK, li2020erd, meagher2016erdHos, 2oddprimes,  spiga2019erdHos}.  In these works, two permutations $g,h \in G  \leq \sym(V)$ are said to be \textsl{intersecting} if $g(v) =h(v)$ for some  element $v \in V$, and the main problem is to determine the largest set of permutations in which any two are intersecting. 
If the subgroup $H \leq G$ is the stabilizer of a point under the action of $G$ on $V$, then this action is equivalent to $G$ acting on the cosets of $H$. Clearly with this action, the subgroup $H$ is an intersecting set. If there are no intersecting sets of cardinality larger than $|H|$, then the group $G$  is said to have the ``Erd\H{o}s-Ko-Rado (EKR) Property". Indeed, such results are often considered to be generalizations of the 
famous Erd\H{o}s-Ko-Rado Theorem.

The most recent work in this area has turned to looking for groups that do not have the EKR-property, this has 
led to trying to measure how far a group can be from having this property.  One way to do this is to use the \textsl{intersection density of a group}. This is a group parameter, introduced in~\cite{li2020erd}, defined for a transitive group $G\leq 
	\sym(n)$, to be the rational number
	\begin{align*}
		\rho(G) := \max\left\{ \frac{|\mathcal{F}| \, n}{|G|} \ :\ \mathcal{F} \subset G \mbox{ is 
		intersecting} \right\}.
	\end{align*} 
	
Since $G$ is transitive, the orbit-stabilizer lemma implies that the stabilizer of a point in $G$ has order $\frac{|G|}{n}$. Since the stabilizer of a point is an intersecting set, $\rho(G) \geq 1$ for any transitive group $G$.  Further, a transitive permutation group has intersection density 1 if and only if it has the EKR property.  In~\cite{HKKMcyclic} the authors initiate a program aimed at obtaining a deeper understanding of the intersection density of transitive permutation groups, with a focus on groups not having the EKR-property. They find many interesting examples using actions with a cyclic point stabilizer. 
	
	In~\cite{BidySym}, the concept of intersection density was extended to vertex-transitive graphs. A graph $X$ is \textsl{vertex transitive} if its automorphism group $\Aut(X)$ acts transitively on the vertex set of $X$.  The \textsl{intersection density} $\rho(X)$ of a 
	vertex-transitive graph $X$ is the largest intersection density among the transitive subgroups of the automorphism group of the graph. Specifically, the intersection density of a graph $X$ is defined to be the rational number
	\begin{align}
		\rho(X) := \max\left\{ \rho(H) \ :\ H\leq \Aut(X), \mbox{ $H$ transitive} \right\}.
	\end{align}

We note here that the intersection density parameter for vertex-transitive graphs only measures the largest possible intersection density of a transitive subgroup; it does not take into account any smaller intersection densities from other transitive subgroups of automorphism. In \cite{KMP}, the intersection density of vertex-transitive graphs was further refined into the \emph{intersection density array}. Given a vertex-transitive graph $X= (V,E)$, the intersection density array of $X$ is the increasing sequence of rational numbers 
\begin{align*}
	\overline{\rho}(X) := [\rho_1,\rho_2,\ldots,\rho_t],
\end{align*}
for some integer $t\geq 0$, such that for any $i\in\{1,2,\ldots,t\}$, 
there exists a transitive subgroup $K\leq \Aut(X)$ such that $\rho_i = \rho(K)$ and for any transitive subgroup $G\leq \Aut(X)$, there exists $i\in \{1,2,\ldots,n\}$ such that $\rho(G) = \rho_i.$ 
This array gives a more robust way of viewing the intersection property of the automorphism group. For example, the Petersen graph has intersection density $2$, whereas its intersection density array is $[1,2]$. Another interesting example is the Tutte-Coxeter graph;  its intersection density is equal to $\frac{3}{2}$ and its intersection density array is $\left[\frac{3}{2}\right]$. That is, every transitive subgroup of the automorphism group of the Tutte-Coxeter graph has intersection density equal to $\frac{3}{2}$. A vertex-transitive graph $X=(V,E)$ exhibiting this property, i.e. $\overline{\rho}(X) = [\rho_1]$, is called \emph{intersection density stable	}.

In this paper, we continue the work in~\cite{BidySym} to determine the intersection density, and if possible, the intersection density array of the \textsl{Kneser graphs}. These are a well-known family of vertex-transitive graphs. For integers $n$ and $k$, with $n \geq 2k$, the Kneser graph $K(n,k)$ has all the $k$-subsets of $\{1,2,\dots,n\}$ as its vertex set and two vertices are adjacent if they are disjoint. For $n >2k$, it is well-known that $\sym(n)$ is the automorphism group of $K(n,k)$ (this is implied by the EKR theorem, see~\cite[Section 7.8]{godsilroyle}). Since $\sym(n)$ is transitive on the $k$-sets of $\{1,2,\dots,n\}$, the graph $K(n,k)$ is vertex transitive. We want to determine the largest intersection density over all subgroups of $\sym(n)$ that are transitive on the $k$-subsets of $\{1, \dots,n\}$ with $n>2k$. 

There has already been much work done to determine the intersection density of $\sym(n)$ with its action on $k$-sets. The most general result is given by Ellis in~\cite{Ellis}, where it is shown that if $n$ is large relative to $k$, then $\sym(n)$ has intersection density 1 under this action. Ellis conjectured the requirement that $n$ be large relative to $k$ is not necessary. Indeed, this has been confirmed for the smallest values of $k$; for $k=2$ this conjecture is proven in~\cite{2setwise}, and for $k=3$, it is proven in~\cite{3setwise}. 

 It is shown in~\cite{2oddprimes} that the alternating group $\alt(n)$ acting on 2-sets has intersection density 1, provided $n\geq 16$. Using \verb*|Sagemath|, it is not hard to verify that this still holds for $6\leq n\leq 15$; but the group $\alt(5)$ acting on the $2$-sets does not have the EKR property, in fact it has intersection density $2$. Further, in~\cite{BidySym} the authors prove that the alternating group acting on $k$-sets with $k =3,4,5$ also has intersection density 1 for $n > 2k$ (they also determine the intersection density of some sporadic groups that are also transitive on the $k$-sets).

To determine the intersection density of the graph $K(n,k)$, it is necessary to determine the largest intersection density over all subgroups that are transitive on the $k$-sets. The next proposition follows from Theorem 14.6.2 of~\cite{godsil2016erdos} and shows that the largest intersection density will be achieved by a minimal (by inclusion) transitive subgroup.

\begin{prop}
  Let $G$ be a transitive group and let $H$ be a transitive subgroup of $G$ (where $H$ has the same action as $G$).
  The intersection density of $G$ is bounded above by the intersection density of $H$.
  \end{prop}

As stated above, the alternating group acting on $3$-sets has intersection density 1 (see~\cite{BidySym}), so we have the following result.
\begin{thm}
Let $n\geq 7$. If $\alt(n)$ is the minimal transitive subgroup of $\sym(n)$ under its action on 3-sets,  then $K(n,3)$ has intersection density equal to $1$.
\end{thm}

In this paper, we will consider values of $n$ for which the alternating group is not the minimal transitive group on the 3-sets. A group that is transitive on the $k$-sets is called $k$-homogeneous, clearly any group that is 3-transitive is also 3-homogeneous. These groups have been classified and we state two results on 3-homogeneous groups that motivate the choice of groups in this work. The first result is taken from~\cite{MR1373655}.

\begin{thm}
	Let $G\leq \sym(n)$ be $3$-transitive. If $G$ is not equal to $\alt(n)$ with $n \geq 5$ nor $\sym(n)$ with $n \geq 7$, then $G$ is one of: $\asl{d}{2}$, $V_{16}.\alt(7)$ (of degree $16$), $\mathieu{11}$ (of degree $12$), $\mathieu{22}, \Aut(\mathieu{22})$, $\mathieu{23}$, $\mathieu{24}$ or 
\[
				\psl{2}{q} \leq G \leq \pgammal{2}{q},
\]			
with degree $q+1$, where $q$ is a prime power.
\label{thm:classification-transitive}
\end{thm}

\begin{thm}[Kantor~\cite{MR306296}]
Suppose that $G$ is $3$-homogeneous on a set of size $n\geq 6$. Then $G$ is $2$-transitive. Moreover, $G$ is $3$-transitive with the exception of 
\[
\psl{2}{q} \leq G \leq \operatorname{P\Sigma L}(2,q),
\]
where $n-1 = q \equiv 3\ (\operatorname{mod} 4)$;
and 
\[
G \in \{ \agl{1}{8},\agammal{1}{8},\agammal{1}{32}  \}.
\]
\label{thm:classification-homogeneous}
\end{thm}

For many values of $n$, the minimal transitive group on the $2$-sets and 3-sets is $\psl{2}{q}$ or contains $\psl{2}{q}$ as a proper subgroup. In this paper we focus on  these groups with their actions on 2 and 3-sets. We note that some of the minimal subgroups in these two theorems are not transitive on the $3$-subsets. In particular $G = \psl{2}{q}$, where $q \equiv 1 \pmod{4}$ in Theorem~\ref{thm:classification-homogeneous} is not 3-homogeneous.

\section{Background Results}

Our approach to this problem is to build a graph for each group that has the property that the cocliques (independent sets) in the graph correspond exactly to the intersecting sets in the group. Then the size of the maximum cocliques can be determined using algebraic techniques. 

Given a group $G$ and a subset $C\subset G\setminus \{1\}$ with the property that $x^{-1} \in C$ whenever $x\in C$, recall that the \emph{Cayley graph} $\operatorname{Cay}(G,C)$ is the graph whose vertex set is $G$ and two group elements $g$ and $h$ are adjacent if and only if $hg^{-1} \in C$. If $C$ has the additional property that $gxg^{-1} \in C$ for all $g\in G$ and $x\in C$, then we say that the Cayley graph $\operatorname{Cay}(G,C)$ is a \emph{normal Cayley graph.
}
For any permutation group $G$, define the \textsl{derangement graph}, denoted by $\Gamma_G$, to have the elements of $G$ as its vertices, and two vertices are adjacent if they are not intersecting. Then the maximum cocliques in $\Gamma_G$ are exactly the maximum intersecting sets in $G$. We can also consider the complement of the derangement graph (denoted $\overline{\Gamma_G}$); clearly the maximum cliques in $\overline{\Gamma_G}$ are maximum intersecting sets. The derangement graph is a Cayley graph with connection-set equal to the set $\der(G)$ of all derangements (i.e., fixed-point-free permutations) of $G$. Further, since the connection-set $\der(G)$ of $\Gamma_{G}$ is a union of conjugacy classes, $\Gamma_{G}$ is a normal Cayley graph. The derangement graph is also a graph in an \textsl{association scheme}, namely the \textsl{conjugacy class association scheme}. Details are given in~\cite[Chapter 3]{godsil2016erdos}.

Since  $\Gamma_{G}$ is a normal Cayley graph, its eigenvalues of can be calculated using the complex irreducible characters of the group $G$. For details see~\cite[Chapter 14]{godsil2016erdos} or~\cite{BidyThesis} as we only state the formula here.

\begin{theorem}
Let $G$ a permutation group. The eigenvalues of $\Gamma_G$ are
\begin{align}\label{eq:evalues}
\lambda_\chi = \frac{1}{\chi(1)} \sum_{g \in \der(G)} \chi(g)
\end{align}
where $\chi$ is taken over all irreducible characters of $G$.
\end{theorem}

A fascinating aspect of this approach is that frequently the eigenvalues of the derangement graph can be used to determine very effective upper bounds on the size of cocliques and cliques in the derangement graphs. The next results are two such bounds, before stating them we need some notation. For a graph $X$ on $n$ vertices, a real symmetric $n \times n$ matrix $M$ is \textsl{compatible} with $X$ if $M_{u,v}=0$ 
whenever $u$ and $v$ are non-adjacent vertices in $X$. 
The \textsl{adjacency matrix} of a graph $X$ on $n$ vertices is an $n \times n$ 01-matrix with $(u,v)$-entry equal to 1 if and only if $u$ and $v$ are adjacent in $X$. The adjacency matrix is denoted by $A(X)$ and is an example of a matrix compatible with $X$.
The sum of all the entries of a matrix $M$ will be denoted by $\elsm(M)$ and the trace by $\tr(M)$. In a graph $X$, the size of the largest coclique is denoted by $\alpha(X)$, and the size of the largest clique by $\omega(X)$. The all ones vector will be denoted by $\mathbf{1}$, and $J$ will represent the all ones matrix (the sizes will be clear from context).

\begin{thm} [Weighted Ratio Bound (Theorem~2.4.2 \cite{godsil2016erdos})]~\label{thm:wtRatio} Let $X$ be a
  connected graph.  Let $A$ be a matrix compatible with $X$ that has
  constant row and column sum $d$.

 If the least eigenvalue of $A$ is $\tau$, then
\[
  \alpha(X) \le\frac{|V(X)|}{1-\frac{d}{\tau}}.
\]
\end{thm}

\begin{theorem}[Theorem~3.7.1 \cite{godsil2016erdos}]\label{thm:LP-clq}
	Let $\mathcal{A} = \{A_0=I, A_1, \dots, A_d\}$ be an association scheme with $d$ classes and let $X$ be a graph with adjacency matrix $A(X) = \sum_{i \in T} A_i$, where $T \subset \{1,2,\ldots,d\}$.
        If $C$ is a clique in $X$, then
	\[
		|C| \le \max_{M \in \mathcal{M} } \frac{\elsm(M)}{\tr(M)}
	\]
where $\mathcal{M}$ is the set of all positive semidefinite matrices in $\mathbb{C}[\mathcal{A}]$ that are compatible with $X$.
\end{theorem} 

The next result gives a simple method to test if a subgroup is an intersecting set.
\begin{lem}\label{lem:intersectingsubgroup}
If $H \leq \sym(n)$ is a subgroup with no derangements, then $H$ is an intersecting set.
\end{lem}
\begin{proof}
If $g, h \in H$, then $g h^{-1} \in H$ and is not a derangement. This means that $gh^{-1}$ has a fixed point, thus $g$ and $h$ are intersecting.
\end{proof}

Finally, we state a well-known result about the transitivity of $\psl{2}{q}$ on the 3-sets that proves it is not $3$-homogeneous.
\begin{lem}
If $q \equiv 1 \pmod{4}$, then $\psl{2}{q}$ has two orbits on the 3-sets from $\{1,\dots, q+1\}$.
\end{lem}

In the next three sections we consider groups containing $\psl{2}{q}$ acting on 3-sets of $\{1,2,\dots, q+1\}$.
The first of these sections considers $\pgl{2}{q}$ where $q$ is even (in this case $\psl{2}{q} =\pgl{2}{q}$). 
Section~\ref{sec:1mod4} considers two groups that contain $\psl{2}{q}$ with $q \equiv 1 \pmod{4}$. In Section~\ref{psl}, we consider $\psl{2}{q}$ for $q \equiv 3 \pmod{4}$, since this is when $\psl{2}{q}$ is transitive.  
Section~\ref{2sets} discusses $\psl{2}{q}$ on the 2-sets. Section~\ref{intrans} briefly considers the intransitive action of $\psl{2}{q}$ on the 3-sets when $q \equiv 1 \pmod{4}$.

%%%%%%%%%%%%%%%%%%%%%%%%%%%%%%%%%%%%%%%%%%%%%%%%
\section{ $\pgl{2}{q}$ acting on $K(q+1,3)$ where $q$ even} 

For $q$ even, we can determine the exact intersection density of $\pgl{2}{q}$ acting on the 3-sets from $\{1,\dots, q+1\}$.

\begin{thm}\label{thm:qeven}
Let $q$ be even. The intersection density of $\pgl{2}{q}$ acting on the $3$-sets is %the following
\begin{enumerate}
\item $\frac{q}{6}$,   if $q = 2^{2 \ell + 1}$;
\item $\frac{q}{2}$,  if $q = 2^{2\ell}$.
\end{enumerate}
\end{thm}

We prove this result in two lemmas. The first is a construction of an intersecting set of the required size.

\begin{lem}
Consider the action of $\pgl{2}{q}$ on the 3-sets from $\{ 1,\dots, q+1\}$.
If $q =2^{2\ell}$ there is an intersecting set of size $3q$, and if $q =2^{2\ell + 1 }$ there is an intersecting set of size $q$. 
\label{lem:max-intersecting-set-PGL-even}
\end{lem}
\begin{proof}
For all $q = 2^n$ the subgroup of $\pgl{2}{q}$ generated by the matrices of the form
\[
\begin{pmatrix} 1 & a \\ 0 & 1\end{pmatrix},
\]
with $a \in \mathbb{F}_q$,  is a subgroup in which all non-identity elements have order 2. Thus, by Lemma~\ref{lem:intersectingsubgroup}, these form an intersecting set under this action.

If $q = 2^{2\ell}$ then there is an $x \in \mathbb{F}_q$ with $x^3=1$.
The set of all matrices of the forms
\[
\begin{pmatrix}
1 & a  \\ 0 & 1 \\
\end{pmatrix},
\quad 
\begin{pmatrix}
x & a  \\ 0 & x^2 \\
\end{pmatrix},
\quad
\begin{pmatrix}
x^2 & a  \\ 0 & x \\
\end{pmatrix}
\]
with $a \in \mathbb{F}_q$,  is a subgroup with size $3q$. Each of these vertices either has order 3, or has order 2 and fixes a point, so each fixes a 3-set. Thus, by Lemma~\ref{lem:intersectingsubgroup}, these matrices form an intersecting set under this action.
\end{proof}

Using Theorem~\ref{thm:LP-clq}, we can show that the sets given in Lemma~\ref{lem:max-intersecting-set-PGL-even} are the largest possible intersecting sets under the action of $\pgl{2}{q}$ on the 3-sets. We note that the stabilizer of $\pgl{2}{q}$ acting on the 3-sets is isomorphic to $\sym(3)$, for any prime power $q$. Henceforth, we denote the stabilizer of a point of $\pgl{2}{q}$ acting on the 3-sets by $H_q$.  Let $X_q$ be the complement of the derangement graph under this action. So the vertices of $X_q$ are the elements of $\pgl{2}{q}$ and two vertices $g, h$ are adjacent if $gh^{-1}$ is conjugate to an element in $H_q$.  A clique in this graph is an intersecting set under the action on the 3-sets. 

The graph $X_q$ is in the conjugacy class association scheme of $\pgl{2}{q}$; we denote this association scheme by $\mathcal{A}$. The matrix in $\mathcal{A}$ that corresponds to the conjugacy class of order 2 elements in $H_q$ will be denoted by $A_1$, and $A_2$ will denote the matrix corresponding to the conjugacy class of the order 3 elements. This means that $A(X_q) = A_1 + A_2$. 
By Theorem~\ref{thm:LP-clq}, any clique in $X_q$ is bounded by the maximum of 
\[
\frac{\elsm(M)}{\tr(M)}
\]
taken over all positive semi-definite matrices $M$ of the form $M = dI + aA_1 + bA_2$.

\begin{lem}\label{lem:even}
Consider the action of $\pgl{2}{q}$ on the 3-sets of $\{1,2,\dots,q+1\}$.
If $q =2^{2\ell}$, then an intersecting set under this action has at most $3q$ elements, and
if $q =2^{2\ell + 1}$ an intersecting set has at most $q$ elements.
\end{lem}
\begin{proof}
 We will apply Theorem~\ref{thm:LP-clq}. Let $\mathcal{A}$ be the conjugacy class association scheme for $\pgl{2}{q}$. We will first find a positive semi-definite matrix $M \in \mathbb{C} [\mathcal{A}]$ that is compatible with the complement of the derangement graph of $\pgl{2}{q}$ under this action, and then we show that $\frac{\elsm(M)}{\tr(M)}$ equals the bounds in the lemma.
 
 Let $X_q$ be the complement of the derangement graph, so $X_q$ is the graph with the elements of $\pgl{2}{q}$ as its vertices and two vertices $g, h$ are adjacent if $gh^{-1}$ is conjugate to an element in $H_q$ (where $H_q$ is the stabilizer of a point under this action). 
  
Let $C_1$ be the conjugacy class of $\pgl{2}{q}$ that contains the elements in $H_q$ of order two, and $C_2$ the conjugacy class that contains the elements of order 3. Define $A_1$ to be the matrix with rows and columns indexed by the element of $\pgl{2}{q}$ and the entry $(g,h)$ equal to 1 if $gh^{-1} \in C_1$ and 0 otherwise; a matrix $A_2$ is defined similarly, but for $C_2$. Both $A_1$ and $A_2$ are matrices in $\mathcal{A}$. The adjacency matrix of $X_q$ is equal to $A_1 + A_2$, and any matrix in $\mathbb{C}[\mathcal{A}] $ compatible with $X_q$ has the form $M=dI + aA_1 +bA_2$. 
If we set $v = | \pgl{2}{q} | $, then
\[
\frac{\elsm(M)}{\tr(M)} =    \frac{v (d + a |C_1| + b |C_2|) }{vd} =  1+ \frac{a}{d} |C_1| + \frac{b}{d} |C_2|.
\]
So we need to find values of $\frac{a}{d}$ and $\frac{b}{d}$ so that the eigenvalues of $M$ are non-negative and $\frac{\elsm(M)}{\tr(M)} $ is maximized.
 
The eigenvalues of $A_1$ and $A_2$ can be calculated easily from the character table of $\pgl{2}{q}$, as the eigenvalue of $A_i$ is simply 
 \[
\lambda_{\chi} (A_i) =  \frac{ |C_i| \chi(c_i) } {\chi(\id)}
\]
where $c_i \in C_i$ and $\chi$ an irreducible character of $\pgl{2}{q}$.

First consider when $q =2^{2\ell}$. The value of all the irreducible characters of $\pgl{2}{q}$ on these two conjugacy classes are known, and recorded in the table below using the notation of~\cite{adams2002character}.

\def\arraystretch{1.3}
\begin{table}[h]
\begin{center}
	\begin{tabular}{|c|cc ccc|} \hline
		Character & $\rho(1)$  & $\rho(\alpha)$  & $\overline{\rho}(1)$  & $ \rho'(1)$ &  $\pi(\chi)$ \\
		    Degree & $q+1$ & $q+1$ & $q$ & $1$  & $q-1$ \\ \hline 
		value on $C_1$ (order 2) &  $1$ & $1$  & $0$  & $1$ & $-1$  \\
		eigenvalue of $A_1$ &  $q-1$ & $q-1$& $0$ &  $q-1$ & $-(q+1)$\\ \hline
		value on $C_2$ (order 3) & $2$ & $-1$ & $1$ & $1$ & $0$  \\
		eigenvalue of $A_2$ & $2q$ &  $-q$ & $q+1$ & $q(q+1)$ & $0$ \\ \hline
	\end{tabular}
\end{center}
\caption{Partial character table for $\pgl{2}{2^{2\ell}}$, with eigenvalues for $A_1$ and $A_2$.}
\end{table}

By Theorem~\ref{thm:LP-clq}, a bound for the size of the cliques in $X_q$ is given by the solution to the following linear program (we use $x$ and $y$ in place of $\frac{a}{d}$ and $\frac{b}{d}$).

\begin{align}
	\begin{split}
	\mathsf{Maximize: } &\ \ \  1 + x (q^2-1) + y q(q+1),\\
	\mathsf{Subject \, to } & 	\\
	& 
	\begin{aligned}
           & -1 \leq x (q-1) +2yq \\
           & -1 \leq x (q-1) - y q \\
           & -1 \leq y(q+1) \\
           & -1 \leq -x(q+1) \\
           %This was -y(q+1)!!!
	\end{aligned}
	\end{split} %\label{general-linear-program}
\end{align}

It is straight-forward to see that this is maximized at $x = \frac{1}{q+1}$ and $y =\frac{2}{q+1}$ to give a maximum value of $3q$.

For $q = 2^{2\ell+1}$, the values that the irreducible characters take on the conjugacy classes with order 2 and 3 are given in the table below.

\def\arraystretch{1.3}
\begin{table}[h]
\begin{center}
	\begin{tabular}{|c|cc ccc|} \hline
		Character & $\rho(\alpha)$  & $\overline{\rho}(1)$  &  $\rho'(1)$ &  $\pi(1)$ &  $\pi(\chi)$ \\
		    Degree & $q+1$ & $q$ & $1$  & $q-1$ & $q-1$ \\ \hline 
		value on $C_1$ (order 2) & $1$  & $0$  & $1$ & $-1$  & $-1$ \\
		eigenvalue of $A_1$ &  $q-1$  & $0$ &  $q-1$ & $-(q+1)$ & $-(q+1)$\\ \hline
		
		value on $C_2$ (order 3) & $0$ & $-1$ & $1$ & $-2$ & $1$  \\
		eigenvalue of $A_2$ & $0$ & $-(q-1)$ & $q(q-1)$ &  $-2q$ & $q$\\ \hline
	\end{tabular}
\end{center}
\caption{Partial character table for $\pgl{2}{2^{2\ell+1}}$, with eigenvalues for $A_1$ and $A_2$.}
\end{table}

The solution of the following linear optimization is a bound on the size of the maximum clique in $X_q$.

\begin{align}
	\begin{split}
	\mathsf{Maximize } &\ \ \  1 + x (q^2-1) + y q(q-1),\\
	\mathsf{Subject \, to } & 	\\
	& 
	\begin{aligned}
           & -1 \leq x (q-1) \\
           & -1 \leq -y (q-1) \\
           & -1 \leq -x (q+1) - 2yq \\
           & -1 \leq -x (q+1)+ yq. \\
	\end{aligned}
	\end{split} %\label{general-linear-program}
\end{align}

It is straight-forward to solve this linear program, it is maximized at $x = \frac{1}{q+1}$ and $y=0$, giving a maximum value of $q$.
\end{proof}

%%%%%%%%%%%%%%%%%%%%%%%%%%%%%%%%%%%%%%%%%%%%%%%%
\section{The subgroups containing $\psl{2}{q}$ on $K(q+1,3)$ when $q \equiv 1 \pmod 4$} 
\label{sec:1mod4}

In this section, we consider the subgroups of the automorphism group of $K(q+1,3)$ containing $\psl{2}{q}$ for a prime power $q \equiv 1 \pmod 4$. In this case, $\psl{2}{q}$ is intransitive, so we consider the two minimally transitive subgroups of the automorphism group of $K(q+1,3)$ containing it. The first one of these groups of course is $\pgl{2}{q}$. The other minimally transitive group is described as follows. If $q = p^k$ for some even number $k$ and an odd prime $p$, then the outer automorphism group of $\psl{2}{q}$ is $\langle \alpha \rangle \times \langle \tau\rangle\cong C_2 \times C_k$, where $\alpha$ and $\tau$ have order $2$ and $k$, respectively. The other minimally transitive subgroup containing $\psl{2}{q}$ is the group $\psl{2}{q} \langle  \alpha \tau^{\frac{k}{2}}\rangle$. For both groups the stabilizer of a 3-set has size 6 and is isomorphic to $\sym(3)$. 

We start this section with a note about the structure of the derangement graph $\Gamma_{G}$ where 
$G$ is either $\pgl{2}{q}$, or  $\psl{2}{q} \langle  \alpha \tau^{\frac{k}{2}}\rangle$.
A graph $X = (V(X),E(X))$ is a \textsl{join} of two vertex-disjoint graphs $Y$ and $Z$, denoted $X = Y \vee Z$, if $V(X) = V(Y) \cup V(Z)$ and $E(X) = E(Y) \cup E(Z) \cup \left\{ (y,z) \mid y \in V(Y), z \in V(Z) \right\}$. That is, $X$ is obtained by taking the disjoint union of $Y$ and $Z$, and adding all the possible edges between the vertices of $Y$ and $Z$. 
We will show that  $\Gamma_{G}$ is a join, to prove this we need to define an operation on the 3-sets.

The group $\psl{2}{q}$ acts on the lines (i.e., $1$-dimensional subspaces of $\mathbb{F}_q^2$). These lines can be represented as homogeneous coordinates of the form $u = (u_1,u_2)$.
For any two vectors $u, v$ we can define  $D(u,v) = u_1 v_2 - u_2 v_1$ and for any ordered triple of lines $(u,v,w)$
consider
\[
D(u,v,w) : = D(u,v) D(v,w) D(w,u).
\]
Note that 
\[
D(v,u,w) = D(v,u) D(u,w) D(w,v) = -D(u,v) D(v,w) D(w,u) = -D(u,v,w).
\]
This value is invariant under scalar multiplication and the action of $\psl{2}{q}$. If $q \equiv 1 \pmod{4}$ then $-1$ 
is a square and for every triple this value will be either a square or a non-square--this shows that the action of $\psl{2}{q}$ has 2 orbits (of equal length by normality) on the 3-sets when $q \equiv 1 \pmod{4}$. Consider $G \in  \{ \pgl{2}{q}, \psl{2}{q} \langle  \alpha \tau^{\frac{k}{2}}\rangle \}$. Any element of $G \backslash \psl{2}{q}$ will map 3-sets from one orbit to the other orbit. This implies that the elements in $G \backslash \psl{2}{q}$ are all derangements, and all the elements of $G$ that fix a 3-set are contained in the subgroup $\psl{2}{q}$. Of these elements, the ones of order 2 are always contained in a single conjugacy class of $\psl{2}{q}$, that we denote by $C_1$; this is also a conjugacy class in $G$. The elements of order 3 are also contained in a single conjugacy class denoted $C_2$, unless $q=3^\ell$. If $q =3^\ell$, then the elements of order 3 are split between two conjugacy classes, $C_2'$ and $C_2''$ in $\psl{2}{q}$ and these conjugacy classes are closed under inversion if and only if $\ell$ is even. For all values of $q$, the elements of order 3 are all contained in a single conjugacy class $C_2$ of $G$. In the case that $q$ is a power of 3, $C_2 = C_2' \cup C_2''$. 

For any $g \in \psl{2}{q}$ and any $h \in G \backslash \psl{2}{q}$, it follows that $g h^{-1} \in G \backslash \psl{2}{q}$. Since all the elements in $G \backslash \psl{2}{q}$ are derangements, this means that the vertices corresponding to $g$ and $h$ are adjacent in $\Gamma_G$. From this we can see $\Gamma_G = X \vee X$ where $X$ is the subgraph induced by permutations $\psl{2}{q}$. In $X$, vertices $g, h$ are adjacent if $gh^{-1}$ is not conjugate, in $G$, to an element in the stabilizer of a 3-set. So vertices $g, h \in X$ are adjacent if $gh^{-1}$ is not in one of the conjugacy classes, $C_1$ or $C_2$ in $G$. 
This shows that $\Gamma_{\pgl{2}{q}}$ and $\Gamma_{\psl{2}{q} \langle  \alpha \tau^{\frac{k}{2}}\rangle }$ are isomorphic since they are both isomorphic to a join of two copies of $X$. Further, since a coclique of $\Gamma_{G} = X \vee X$ must lie in a copy of $X$, the cocliques in $\Gamma_{G}$ are exactly the cocliques in the subgraph induced by the elements of $\psl{2}{q}$.

\begin{lem}
Let $q \equiv 1 \pmod{4}$ and $G \in  \{ \pgl{2}{q},  \psl{2}{q} \langle  \alpha \tau^{\frac{k}{2}}\rangle\}$. Then $\Gamma_{G}$, with the action on the 3-sets, is the join of two copies of a graph $X$.  The vertices of $X$ are the elements of $\psl{2}{q}$ and two elements are adjacent if they are not conjugate, in $G$, to an element in the stabilizer of a 3-set.
\end{lem}

Note that if $q = 3^{2\ell}$, then $X$ is not the same as the derangement graph from the action of $\psl{2}{q}$ acting on the 3-sets (which is the action considered in the next section). In the graph $\Gamma_{\psl{2}{q}}$, with the action on the 3-sets, two vertices $g,h$ are adjacent if $gh^{-1}$ is in $C_1$ or $C_2'$ (but not in $C_2''$). We will consider the intersection density of the intransitive subgroup $\psl{2}{q}$ in Section~\ref{intrans}.

For the remainder of this section, we only consider $\pgl{2}{q}$, since the derangement graphs are isomorphic, the same results will hold for $\psl{2}{q} \langle  \alpha \tau^{\frac{k}{2}}\rangle$. The size of $\pgl{2}{q}$ is $(q-1)q(q+1)$, and the stabilizer of a 3-set is $H_q \cong \sym(3)$. We start with the case where $q \equiv 1 \pmod 4$ is a power of $3$. We note that this implies that $q = 3^{2 \ell}$. 

\begin{thm}\label{thm:powerof3}
If  $q =3^{2\ell}$, for some positive integer $\ell$, then the intersection density of $\pgl{2}{q}$ with its action on 3-sets of $\{1, \dots, q+1\}$ is $\frac{q}{3}$.
\end{thm}
\begin{proof}
First, the set of matrices of the forms
\[
\begin{pmatrix} 1 & x \\ 0 & 1\end{pmatrix},
\quad
\begin{pmatrix} 1 & x \\ 0 & -1 \end{pmatrix},
\]
with $x \in \fld_q$,  form a subgroup that is an intersecting set of size $2q$ with the action of $\pgl{2}{q}$ on the 3-sets. 

Next we will use Theorem~\ref{thm:LP-clq}, with the method in the previous section to show that $2q$ is an upper bound on the size of a clique in $X_q$. 
Again, define a graph $X_q$ whose vertices are the elements of $\pgl{2}{q}$ and two vertices $g, h$ are adjacent if $gh^{-1}$ is in one the conjugacy classes that contain elements from $H_q$.  Using the same notation as in the proof of Lemma~\ref{lem:even}, this implies that $A(X_q) = A_1 + A_2$.

Table~\ref{tab:pgl32l} gives the values of all the irreducible characters of $\pgl{2}{q}$, this is taken from~\cite{adams2002character} and we use the same notation.
\def\arraystretch{1.3}
\begin{table}
\begin{center}
\begin{tabular}{|c| cc cc c|} \hline
Character & $\rho(\alpha)$  & $\rho(\alpha)$  & $\overline{\rho}(\alpha)$  & $\rho'(\alpha)$  & $\pi(\chi)$ \\
   Dimension & $q+1$ & $q+1$ & $q$ & $1$  & $q-1$ \\ \hline 
value on $C_1$ &  $2$ & $-2$  & $1$  & $1$ &  $0$ \\
eigenvalue of $A_1$ & $ q$ & $-q$ & $\frac{q+1}{2}$  & $\frac{q(q+1)}{2} $ & $0$ \\ \hline
value on $C_2$ & $1$ & $1$ & $0$ & $1$  & $-1$ \\
eigenvalue $A_2$ & $q-1$ &  $q-1$ & $0$ & $q^2-1$ & $-(q+1)$ \\ \hline
\end{tabular}
\end{center}
\caption{Partial character table for $\pgl{2}{3^{2\ell} }$, with eigenvalues for $A_1$ and $A_2$.}\label{tab:pgl32l}
\end{table}

The linear optimization problem we need to solve is the following:
\bigskip

\begin{align}
	\begin{split}
	\mathsf{Maximize } &\ \ \  1 +  x (q^2-1) + y \frac{q(q+1)}{2},\\
	\mathsf{Subject \, to } & 	\\
	& 
	\begin{aligned}
	   & -1 \leq x \frac{(q^2-1)}{q+1} + 2 y \frac{q(q+1)}{2(q+1)}  = x (q-1) + y q  \\
           & -1 \leq x \frac{(q^2-1)}{q+1} - 2 y \frac{q(q+1)}{2(q+1)} = x (q-1) - y q  \\           
           & -1\leq y \frac{q(q+1)}{2q}  = y \frac{q+1}{2}  \\
           & -1 \leq x (q^2-1) + y \frac{q(q+1)}{2}    \\
           & -1 \leq x  \frac{ -(q^2-1)}{q-1}  = -x(q+1)  \\
	\end{aligned}
	\end{split} % \label{general-linear-program}
\end{align}

This linear program can be easily solved to see that the objective function is maximized at $x = \frac{1}{q+1}$ and $y =\frac{2}{q+1}$ with a value of 
\[
1 + \frac{1}{q+1} (q^2-1) + \frac{2}{q+1} \frac{q(q+1)}{2} = 1 + (q-1) + q = 2q.\qedhere
\]
\end{proof}

%%%%%%%%%%%%%%%%%%%%%%%%%%%%%%%%
Next, we show that there is an intersecting set of size twice the order of the stabilizer of a 3-set in $\pgl{2}{q}$, for any $q \equiv 1 \pmod 4$.
\begin{lem}\label{lem:alt}
If $q \equiv 1 \pmod{4}$, then there is an intersecting set of size $12$ in $\pgl{2}{q}$ with the action on the 3-sets of $\{1, \dots ,q+1\}$.
\end{lem}
\begin{proof}
The group $\pgl{2}{q}$ has a subgroup isomorphic to $\alt(4)$.  Each of the elements in $\alt(4)$ have order either 2 or 3, so all the element of order 3 fix at least one 3-set. Further, since $q \equiv 1 \pmod{4}$, any element with order 2 has $\frac{q-1}{2}$ 2-cycles. This means any such permutation has two fixed points, so it will also fix a 3-set. Thus by Lemma~\ref{lem:intersectingsubgroup}, this subgroup is an intersecting set.
\end{proof}

We conjecture that the intersecting sets in Lemma~\ref{lem:alt} are the largest.

\begin{conj}
If $q$ is not a power of $2$ or $3$, and $q \equiv 1 \pmod{4}$, then the intersection density of $\pgl{2}{q}$ acting on the 3-sets from $\{1, \dots, q+1\}$ is $2$.
\end{conj}

%%%%%%%%%%%%%%%%%%%%%%%%%%%%%%%%%%%%%%
\section{$\psl{2}{q}$ on $K(q+1,3)$ when $q \equiv 3 \pmod 4$} 
\label{psl}

Next we will consider the group $\psl{2}{q}$ acting on 3-sets from $\{1,\dots, q+1\}$ where $q \equiv 3 \pmod 4$. This action is transitive, and the stabilizer of a 3-set in $\psl{2}{q}$ has size $3$ and is isomorphic to $\mathbb{Z}_3$.

First we consider when $q$ is a power of $3$. Since we are only considering $q \equiv 3 \pmod{4}$ in this section, $q$ must be an odd power of $3$.

\begin{thm}\label{thm:oddpowerof3}
If $q = 3^{2\ell+1}$ the intersection density of $\psl{2}{q}$ under the action on $3$-sets is $\frac{q}{3}$.
\end{thm}
\begin{proof}
First there is an intersecting set of size $q$ given by the subgroup of matrices with the form
\[
\begin{pmatrix} 1 & x \\ 0 & 1\end{pmatrix}
\]
where $x \in \fld_q$. 

We define the derangement graph $\Gamma_{\psl{2}{q}}$, as usual, with the elements of $\psl{2}{q}$ as its vertices, and two vertices $g, h$ are adjacent if $gh^{-1}$ is a derangement. A maximum intersecting set of the group is a maximum coclique in this graph. Using Theorem~\ref{thm:wtRatio}, we will show a coclique in $\Gamma_{\psl{2}{q}}$ is no larger than $q$, so the subgroup above is the largest possible intersecting set.
 In~\cite{adams2002character} the conjugacy classes of $\psl{2}{q}$ are grouped into families; the families denoted by $c_3(x)$ (with $x \neq 1$ and $x^2 \neq -1$) and $c_4(x)$ are exactly the derangements under this action. Below, we record the sums of the values of the irreducible characters over the conjugacy classes in each family.

\begin{table}[h]
\def\arraystretch{1.3}
\begin{center}
\begin{tabular}{|c |ccc ccc|} \hline
Character &    $\rho'(1)$  & $\overline{\rho}(1)$ &  $\rho(\alpha)$   &  $\pi(\chi)$& $\pi(\chi)$  & $\omega_0^{\pm}$  \\
Dimension              &   1 & $q$ & $q+1$ & $q-1$ & $q-1$ & $\frac{q-1}{2}$ \\ \hline
value on $c_3 (x)$ &  $\frac{q-3}{4}$  & $\frac{q-3}{4}$ & $-1$ &  0 & 0 & 0 \\ 
eigenvalue &  $\frac{q(q+1)(q-3)}{4}$  & $\frac{(q+1)(q-3)}{4}$ & $-q$ & 0 & 0 & 0 \\ \hline

value on $c_4(x)$ &  $\frac{q-3}{4}$  &  $-\frac{q-3}{4}$ & 0 &  0 & 2 &  0 \\ 
eigenvalue &        $\frac{q(q-1)(q-3)}{4}$  &  $-\frac{(q-1)(q-3)}{4}$ & 0 & 0 & $2q$ & 0 \\ \hline 
\end{tabular}
\end{center}
\caption{Partial character table $\psl{2}{3^{2\ell+1}}$.}
\end{table}

 From this, it is straight-forward to find the eigenvalues of a matrix compatible with $\Gamma_{\psl{2}{q}}$.
If we set the weight on the conjugacy classes of type $c_3(x)$ to be $a= \frac{1}{q}$ and weight of the conjugacy class of type $c_4(z)$ to be $b  = \frac{(q+3) }{q(q-3) } $, then the largest eigenvalue of the weighted adjacency matrix is
\[
\frac{1}{q}  \frac{q(q+1)(q-3)}{4}  +   \frac{(q+3) }{q(q-3) }  \frac{q(q-1)(q-3)}{4} 
=\frac{q^2-1}{2} - 1 ,
\]
and the smallest is $-1$ (from both $\overline{\rho}(1)$ and $\rho(\alpha)$).
Then the ratio bound gives the size of the largest coclique
\[
\frac{(q-1)q(q+1)}{2(\frac{q^2-1}{2} )} = q.\qedhere
\]
\end{proof}

This group action has also been considered in~\cite{HKKMcyclic}, where an exact value has been determined for some values of $q$.

%M=dI + aA_1 +bA_2$

\begin{thm}[Theorem 6.1, \cite{HKKMcyclic}] \label{thm:q1mod3}
Consider $\psl{2}{q}$ with its action on 3-sets. If $q = p^\ell$ with $q \equiv 1 \pmod{3}$ then
\[
\rho(G) = 
\begin{cases}
4/3 & \textrm{ if } p\neq5, \\
2 & \textrm{ if } p = 5. \\
\end{cases}
\]
\end{thm}

\medskip

If $q$ is such that $q \equiv 3 \pmod 4$ and $q \equiv 2 \pmod{3}$ then $q \equiv 11 \pmod{12}$.
For $q \equiv 11 \pmod{12}$, a simple calculation shows that $q^2 \equiv 1 \pmod{5}$ or $q^2 \equiv 4 \pmod{5}$.

\begin{lem}
If $q^2 \equiv 1 \pmod{5}$ then the density of $\psl{2}{q}$ with its action on 3-sets is at least $4/3$. 
\end{lem}
\begin{proof}
If $q^2 \equiv 1 \pmod{5}$ then $\psl{2}{q}$ contains a copy of $\alt(5)$. The subgroup $\alt(5)$ has an intersecting set of size 4.  An example of such a set is 
\[
\{id, (1,2,3), (1,2,4), (1,2,5) \}.
\]
The subgraph induced by this subset of $\psl{2}{q}$ is a clique of size 4 under this action.
\end{proof}

We end this section with a conjecture on the intersection density of the group $\psl{2}{q}$ with its action on the 3-sets, when $q^2 \equiv \pm 1\pmod{5}$.

\begin{conj}
If $q^2 \equiv 1 \pmod{5}$ then the intersection density of $\psl{2}{q}$ with its action on 3-sets is $4/3$; if $q^2 \equiv 4 \pmod{5}$ then the intersection density of this action is $1$.
\end{conj}

%%%%%%%%%%%%%%%%%%%%%%%%%%%%%%%%%%%%%%
%%%%%%%%%%%%%%%%%%%%%%%%%%%%%%%%%%%%%%

\section{$\psl{2}{q}$ acting on the Kneser graph $K(q+1,2)$}
\label{2sets}

Both the groups $\pgl{2}{q}$ and $\psl{2}{q}$ are transitive subgroups of the automorphism group of $K(q+1,2)$ where $q$ is a prime power. We will only consider $\psl{2}{q}$ since it would have the larger intersection density of the two.
If $q$ is odd, the size of the stabilizer of a point on $\psl{2}{q}$ under the action on 2-sets has size
\[
\frac{(q-1)q(q+1)}{2}  \binom{q+1}{2} ^{-1} = q-1,
\]
and if $q$ is even it is
\[
(q-1)q(q+1)  \binom{q+1}{2} ^{-1} = 2(q-1).
\]

\begin{lem}
For $q$ even, the group $\psl{2}{q}$ acting on 2-sets from $\{1,\dots,q+1\}$ has intersection density $\frac{q}{2}$.
\end{lem}
\begin{proof}
Since $q$ is even, any element that fixes exactly one point has order 2, and also fixes a 2-set. 
Clearly any element with two fixed points also fixes a 2-set. Thus every element in the stabilizer of a point in $\psl{2}{q}$, in its action on the points $\{1,\dots , q+1\}$, also fixes a 2-set from $\{1,\dots , q+1\}$. So the stabilizer of a point in the natural action is also an intersecting set under the action on 2-sets and it has size $q(q-1)$.

We will use the ratio bound, Theorem~\ref{thm:wtRatio}, to show that this is the largest possible intersecting set. Still using the notation of~\cite{adams2002character}, only the conjugacy classes of type $c_4(z)$ are derangements with this action and the eigenvalues of the derangement graph are
\[
\left\{ \frac{q^2(q-1)}{2}, \quad 0, \quad -\frac{q(q-1)}{2}, \quad q \right \}.
\]
The ratio between the largest and the smallest eigenvalue is $-q$, so by the ratio bound, Theorem~\ref{thm:wtRatio}, the size of a maximum coclique cannot be any larger than
\[
\frac{(q-1)q(q+1) }{1 -  (-q) }
=(q-1)q .
\]

This implies the intersection density is
\[
\frac{q^2-q}{ 2(q-1) } = \frac{q}{2}.\qedhere
\]
\end{proof}

For values of $q$ smaller than $32$ we have done some calculations of the intersection density, and based on these calculations we make the following conjecture.

\begin{conj}
For $q \equiv 1 \pmod{4}$ the action of the group $\psl{2}{q}$ on 2-sets has intersection density 1.
\end{conj}

The case for $q \equiv 3 \pmod{4}$ seems more complicated, we did find groups that had larger intersecting sets, but we did not find a general construction.

\begin{lem}\label{lem:counters}
For $q = 7$ the group $\alt(4)$ is an intersecting subgroup in $\psl{2}{q}$ and the intersection density is at least $2$.
For $q = 31$ the group $\alt(5)$ is an intersecting subgroup in $\psl{2}{q}$ and the intersection density is at least $2$.
\end{lem}
The second part of Lemma~\ref{lem:counters} is Example 2.2 in~\cite{li2020erd}, in this paper the authors show that there are no larger intersecting subgroups, but there may be larger intersecting subsets. There are larger values of $q$ for which $\psl{2}{q}$ contains $\alt(4)$ or $\alt(5)$ as intersecting subgroup, but in these groups, the stabilizer under the action on 2-sets is larger than either $\alt(4)$ or $\alt(5)$; so these subgroups do not give an intersection density above $1$. This leaves us with the following question.

\begin{quest}
For $q \equiv 3 \pmod{4}$ and $q>31$, does the group $\psl{2}{q}$ have intersection density larger than 1 with its action on the 2-sets?
\end{quest}

Finally, we also make a conjecture about the intersection density of $\pgl{2}{q}$, even though the maximum intersection density of the Kneser graph would be given by the group $\psl{2}{q}$. From our computations, this conjecture seems true for $q\leq 27$.

\begin{conj}
For $q$ odd the group $\pgl{2}{q}$ with its action on the 2-sets of $\{1,\dots ,q+1\}$ has intersection density 1.
\end{conj}

\section{The intransitive action of $\psl{2}{q}$ on $K(q+1,3)$}
\label{intrans}

If $q \equiv 1 \pmod{4}$, then $\psl{2}{q}$ is not transitive on the 3-sets, but we can consider the action of the $\psl{2}{q}$ on one of its orbits of 3-sets. The stabilizer of a 3-sets under this action is a subgroup $H_q$, isomorphic to $\sym(3)$.
If $q = 3^\ell$ then the conjugacy class of elements of order three splits into 2 conjugacy classes in $\psl{2}{q}$. If $\ell$ is even then these classes are closed under inverses, so all the order three elements of $H_q$ belong to only one of these conjugacy classes; if $\ell$ is odd these conjugacy classes are not closed under inverses. We have done some computer searches for intersecting sets under this action. Our results are recorded in Table~\ref{tab:intrans}.

\begin{table}[h]
\begin{center}
\begin{tabular}{|c | c| c  |c|} \hline
$q$ & Max. intersecting & Lower bound on Density & Intersection density  \\ 
& set found & of  $\psl{2}{q}$ on $H_q$ &  of $\pgl{2}{q}$ on $H_q$ \\ \hline
     5  & 12 & 2 & 2\\
     9 &  15 & 5/2 & 3 \\
    13 & 12 & 2 & 2\\
    17 & 12 & 2 & 2\\
    25 &12 & 2 & 2\\ \hline
\end{tabular}
\end{center}
\caption{Calculations for intersection density of $\psl{2}{q}$ acting on a single orbit of 3-sets when $q \equiv 1 \pmod{4}$.}
\label{tab:intrans}
\end{table}

We consider the group $\pgl{2}{9}$ more carefully; it is an example where there are 2 conjugacy classes of elements of order three and each such class is closed under inverses.

\begin{example}
The group $\psl{2}{9}$ acting on $\sym(3)$ has an intersecting set of size 15.
The elements in this set are the following:

\begin{center}
\begin{tabular}{cc}
 $id$ & \\
$(1,2)(5,10)(6,9)(7,8)$  & $(1,2)(3,4)(5,7)(8,10)$ \\
$(1,10)(2,7)(3,6)(5,8)$  & $(1,7)(2,10)(4,9)(5,8) $ \\
$(1,5)(2,8)(4,6)(7,10)$ & $(1,8)(2,5)(3,9)(7,10) $  \\
$( 1,4,2)(5,8,6)(7,10,9)$ & $(1,2,4)(5,6,8)(7,9,10)$ \\
$(1,6,2)(3,10,7)(4,5,8)$ & $(1,2,6)(3,7,10)(4,8,5)$ \\
$(1,3,2)(5,9,8)(6,10,7)$ & $(1,2,3)(5,8,9)(6,7,10)$ \\
$(1,9,2)(3,8,5)(4,7,10)$  & $(1,2,9)(3,5,8)(4,10,7) $\\
\end{tabular}
\end{center}

Any pair of order 2 elements in a single row generate a subgroup isomorphic to $C_2 \times C_2$. Any other pair of order 2 elements generates a subgroup isomorphic to $\sym(3)$.
Any pair of order 3 elements in a single row are inverses, so generate a copy of $C_3$; any pair of order 3 element not in a row generate a subgroup isomorphic to $\alt(4)$. Finally, an element of order 2 and an element of order 3 generate subgroups isomorphic to either $\sym(3)$ or $\alt(4)$. If an element of order two generates a subgroup isomorphic to $\sym(3)$ with an element of order 3, then the other element of order two in the same row generates a subgroup isomorphic to $\alt(4)$ with the element of order 3.

\end{example}

%%%%%%%%%%%%%%%%%%%%
\section{Conclusions and Further Work}

In this paper, we initiated the study of intersection density of the Kneser graphs $K(n,3)$ and $K(n,2)$. Our main focus is the $3$-homogeneous groups containing the almost simple group $\psl{2}{q}$, for some prime power $q$. The group $\psl{2}{q}$ is not transitive in its action on the $3$-sets when $q \equiv 1 \pmod 4$; for other values of $q$ modulo $4$, it is transitive on $3$-sets.

From Theorems~\ref{thm:classification-transitive} and~\ref{thm:classification-homogeneous}, in the case $n = 2^\ell+1$ for some $\ell$, the minimal 3-transitive subgroups on $\{1, \dots, n\}$ are $\pgl{2}{2^{\ell} }$. Theorem~\ref{thm:qeven} gives the exact values of the intersection density of $K(n,3)$ in this case.

\begin{theorem}
The intersection density of $K(2^\ell+1,3)$ is $3(2^{\ell})$ if $\ell$ is even, and $2^{\ell}$ if $\ell$ is odd.
\end{theorem}

Further, provided that $3^\ell+1 \neq 2^d$ for some $d$, the smallest 3-transitive subgroup is $\pgl{2}{3^{\ell} } $, this, with Theorem~\ref{thm:powerof3} and Theorem~\ref{thm:oddpowerof3}, gives the intersection density of $K(3^{\ell} + 1 , 3)$. In 2004, Mihăilescu \cite{Catalan} proved the Catalan conjecture, which asserts that the only solution to the equation $x^a = y^b+1$, where $(x,y) \in \mathbb{N} \times \mathbb{N}$ and $a,b>1$, is the pair $(x,y)= (3,2)$ with $a =2$ and $b =3$. Therefore, we can see that $3^\ell+1 \neq 2^d$ unless $\ell =1$ and $d=2$.

\begin{theorem}
If $\ell >1$, then the intersection density of $K(3^{\ell} + 1 , 3)$
 is $3^{\ell-1}$.
\end{theorem}

The next project is to calculate the intersection density of $\asl{d}{2}$ with its action on the 3-sets. In particular, we would like to determine if there exists a weighted adjacency matrix for the action of $\asl{d}{2}$ on the 3-sets, for all values of $d$, for which the ratio bound can be used to prove that the stabilizers are the largest intersecting sets.
This means the question of the intersection density of $K(2^\ell,3)$ is still open, as is the question which groups achieve the largest density. 

\begin{quest}
Does the subgroup $\pgl{2}{q}$ give the maximum intersection density of $K(q+1,3)$ for $q$ a prime power with $q \equiv 1 \pmod{4}$?
\end{quest}

In the case that $q \equiv 3 \pmod{4}$, we conjecture that the group $\psl{2}{q}$ gives the maximum intersection density.

\begin{conj} 
For $q$ a prime power with $q \equiv 3 \pmod{4}$, the group $\psl{2}{q}$ gives the maximum intersection density among all subgroups of automorphism of $K(q+1,3)$. If $q^2 = 1 \pmod{5}$ then the intersection density of $K(q+1,3)$ is $4/3$, and if $q^2 = 4 \pmod{5}$ then the intersection density of $K(q+1,3)$ is $1$.
\end{conj}

We would like to determine the intersection density of the other groups acting transitively on the 3-sets. Using \verb|Sagemath|~\cite{sagemath}, we verified that $\agl{1}{8}$ and $\agammal{1}{8}$ with their actions on the $3$-sets both have intersection density 1. The groups $V_{16}.\alt(7)$ has intersection density 1 (via \verb|Sagemath|~\cite{sagemath} and \verb|Gurobi| \cite{gurobi}) and $\agammal{1}{32}$ has intersection density 1, as it is regular. The group $\mathieu{11}(12)$ has intersection density 1 (this was verified by determining that the maximum coclique in the subgraph of the derangement graph induced by the non-neighbours of the identity has size $35$.)

Using \verb| Gurobi |~\cite{gurobi}, there exists a weighted adjacency matrix for the derangement graph so that the ratio bound holds with equality for the groups (with the action on the 3-sets) $\asl{3}{2}$, $\asl{4}{2}$ and $M_{24}$. It can be further determined that the derangement graph for groups $M_{11}$, $M_{22}$ and $M_{23}$ do not have such a weighted adjacency matrix.  We still have yet to determine the intersection density of the Mathieu groups, $M_{22}$ and $M_{23}$, and the group $\Aut(\mathieu{22})$ with their action on the 3-sets. Putting these calculations together, we get the following lemma.

\begin{lem}
	The groups $\agl{1}{8}$, $\agammal{1}{8}$, $\agammal{1}{32}$, $\asl{3}{2}$, $\asl{4}{2}$, $V_{16}.\alt(7)$, $\mathieu{11} (12)$ and $M_{24}$ acting on the 3-sets have intersection density 1.
\end{lem}

As stated in the introduction, the notion of intersection density has been generalized by 
Kutnar, Maru{\v{s}}i{\v{c}} and Pujol in~\cite{KMP} to an intersection density array.
For a vertex-transitive graph, this array consists of the intersection densities for all transitive subgroups of the automorphism group of the graph. The intersection density of a graph is the largest entry in this array, but 
the entire intersection density array for the Kneser graphs is unknown in general. If $\Alt(n)$ is the smallest transitive automorphism of $K(n,k)$, then the intersection density array of $K(n,k)$ is $[1]$; in this case $K(n,k)$ is intersection density stable and has the EKR-property. In this paper we give several examples of Kneser graphs that have intersection density greater than $1$ and hence do not have the EKR property. For these graphs we would like to determine the entire intersection density array. For small Kneser graphs, using  \texttt{GAP}~\cite{GAP4}, we were able to compute the entire array, these are recorded in Table~\ref{ArrayTable}. In the examples we checked, with the exception of $q=8$, all the groups $G$ with $\psl{2}{q} \leq G \leq \pgammal{2}{q}$, the size of the largest coclique in $G$ is the same as the size of largest coclique in $\psl{2}{q}$.

\begin{table}[h]
\begin{center}
\begin{tabular}{|c|c|} \hline
Kneser Graph & Array \\ \hline \hline
%K(4,3) &     [1] \\  \hline %   check alt(4)
%K(5,3) &     [2,1] \\  \hline %k checked
%K(6,3) &     [2,1] \\  \hline %k checked
$K(7,3)$ &     $[1]$ \\  \hline %  check psl(3,2)
$K(8,3)$ &     $[1,4/3]$ \\  \hline %   check psl(3,2) 
$K(9,3)$   &    $[1 ,4/3 ]$ \\ \hline %k checked
$K(10,3)$  &   $[1,3]$  \\  \hline % kchecked
\hline 
%K(12,3)  &   [4/3,1]  \\  \hline % kchecked

$K(17,3)$  & $[1,2,4,8]$ \\ \hline % k checked
$K(26,3)$ & $[1, 2]$  \\ \hline % kchecked
$K(28,3)$ &  $[1,3,9]$ \\ \hline
\end{tabular}
\end{center}
\caption{Intersection Array for some Kneser Graphs.\label{ArrayTable}}
\end{table}

\section{Acknowledgements}
We wish to thank the anonymous referee for their valuable comments which greatly improved this paper.
This research was done while the second author was a Ph.D. student under the supervision of Dr. Karen Meagher and Dr. Shaun Fallat at the Department of Mathematics and Statistics, University of Regina.

\end{document}